\documentstyle[12pt]{article}

\title{\huge\bf Discrete analogues of the Laguerre Inequality}

\author{   {\bf Ilia Krasikov}\\
	   Brunel University\\
	   Department of Mathematical Sciences\\
	   Uxbridge UB8 3PH United Kingdom\\
	   e-mail: mastiik@brunel.ac.uk
       }

\date{}

\newcommand{\RP}{\mbox{$\cal R \cal P \;$}}

\newcommand{\QED}{\hfill$\Box$}

\newcommand{\LP}{\mbox{$\cal L-\cal P \; $}}
\newcommand{\LPO}{\mbox{$\cal L-\cal P (O)\; $}}

\newcommand{\be}{\begin{equation}}
\newcommand{\ee}{\end{equation}}
\newcommand{\proof}{\noindent {\bf Proof.\ \ }}
\newtheorem{lemma}{Lemma}
\newtheorem{theorem}{Theorem}

\newtheorem{conjecture}{Conjecture}

\mathchardef\inn="3232
\renewcommand{\in}{\mbox{$\,\inn\,$}}

\begin{document}


\maketitle
\vspace*{4ex}

\begin{center}
{\bf Abstract}
\end{center}

\begin{quote}
\baselineskip3ex
It is shown that 
$\sum_{j=-m}^m (-1)^j \frac{f(x-j)(f(x+j)}{(m-j)! (m+j)!} \ge 0,$
$m=0,1,...,$
where $f(x)$ is either a real polynomial with only real zeros
or an allied entire function of a special type,
provided the distance between two consecutive zeros of $f(x)$
is at least $\sqrt{4-\frac{6}{m+2}}.$
These inequalities are a surprisingly similar discrete analogue
of higher degree generalizations of
the Laguerre and Turan inequalities. 
Being applied to the classical discrete orthogonal polynomials,
they yield sharp, explicit bounds
uniform in all parameters involved, 
on the polynomials and their extreme zeros.  
We will illustrate it for the case of Krawtchouk polynomials.
\end{quote}
\vspace{16ex}

\noindent


\thispagestyle{empty}
\addtocounter{page}{-1}

\newpage
\baselineskip3.3ex


\renewcommand{\thesection}{\arabic{section}.}
\section{Introduction}
\label{sec1}
A real entire function $\phi (x)$ is in the Laguerre-Polya class \LP 
if it has a representation of the form
$$\phi (x)=c x^m e^{-\alpha x^2+\beta x} 
\prod_{k=1}^{\omega} (1+\frac{x}{x_k}) e^{-x/{x_k}} \, \, \, \, ( \omega
\le \infty ),$$
where $c, \beta, x_k$ are real, $\alpha \ge 0$, $m$ is a nonnegative integer and
$\sum x_k^{-2} < \infty.$\\
We denote by \LPO a subclass of \LP corresponding to $\alpha=0,$
and by $\RP$ a set of real polynomials having only real zeros.
Obviously \RP $\subset$ \LPO $\subset$ \LP.
The well-known Laguerre inequality \cite{lagu} states that
$f'^2-f f'' \ge 0,$
for any $f \in \LP.$
The following remarkable high degree generalization 
\be
\label{pjen}
L_{m}(f(x))=\sum_{j=-m}^{m} (-1)^{m+j} \frac{ f^{(m-j)}(x) f^{(m+j)}(x)}{(m-j)! (m+j)!} \ge 0
\ee
yielding the Laguerre inequality for $m=1$, seems much less known \cite{jensen,patrick1,patrick2}.
M. Patrick \cite{patrick2} used (\ref{pjen}) to obtain 
Turan-type inequalities, which have essentially the same form 
\be
\label{patrick}
\sum_{j=-m}^m (-1)^j \frac{p_{k-j}(x)p_{k+j}(x)}{( m-j)!(m+j)!} \ge 0
\ee
and are valid
for polynomials and entire functions having a
generating function of the Laguerre-Polya class.
Some further generalizations and allied inequalities can be found in 
\cite{dilch, fk3}. A case of entire functions in the Laguerre-Polya class
and a connection between Laguerre and Turan
type inequalities are considered in \cite{craven,csordas,patrick2}.
The aim of this paper is to establish a difference analogue of (\ref{pjen}) and (\ref{patrick})
which is surprisingly similar.\\
Let $x_1 \le x_2 \le ...,$ be the zeros of
$f \in $\LPO, a $mesh$ $M(f)$ of $f$ is defined by
$M(f)=\inf_{i \ge 2}{(x_i-x_{i-1})}.$
Our main result is the following theorem (the case $m=1$ has been given in \cite{k3}).
\begin{theorem}
\label{disin4} Let $f \in$ \LPO and
$M(f) \ge \sqrt{4-\frac{6}{m+2}}, \;$  then
\be
\label{glform}
V_m(f)=\sum_{j=-m}^m (-1)^j \frac{f(x-j)f(x+j)}{( m-j)!(m+j)!} \ge 0
\ee
\end{theorem}
Laguerre and Turan type inequalities have many important applications.
In particular they can be used to find sharp explicit bounds, uniform in all parameters
involved for the classical orthogonal polynomials on the real axis and their extreme zeros.  
The simplest possibility is to convert (\ref{pjen}) into
an algebraic inequality by expressing all the derivatives involved in terms of $f(x)$ and $f'(x),$
using the corresponding second-order differential equation.
In \cite{fk1,fk2,k2,k3, kbes} this type of arguments was applied to the case of Hermite and Krawtchouk
polynomials and the Bessel function. Roughly speaking, inequalities with $m=1$
yield the first term of the corresponding asymptotics, whenever for $m=2$ they
give the correct order for the second term with a slightly weaker constant, probably the best possible
precision achievable by this method. 
For example, for the largest zero $z_k$ of Hermite polynomial $H_k(x)$,
(\ref{pjen}) it gives $z_k \le \sqrt{2k-2}$  with $m=1$ and 
$z_k \le \sqrt{(4k-3k^{1/3}-1)/2} =\sqrt{2k}-\frac{3}{4 \sqrt{2}} k^{-1/6}+O(k^{-1/2}),$ with $m=2.$
Here the constant $\frac{3}{4 \sqrt{2}}$ is about 3 times smaller than the actual asymptotic value.
As well the same approach yields two-sides bounds for the envelop of $|H_k(x)|$ in the oscillatory region 
with the precision $O(k^{-2})$ \cite{fk2}.
We would like to stress
that the method produces inequalities rather than asymptotics
usually obtained by a standard approach. 
Our search for
higher degree discrete analogues of (\ref{pjen}) was motivated by the fact that
the classical discrete orthogonal polynomials satisfy a second order
difference equation of the form
\be
\label{rec1}
p_k(x+1)=b_k(x)p(x)-c_k(x)p_k(x-1)
\ee
with $c_k(x) >0$ for $x$ belonging to the interval of orthogonality \cite{niki85}.
Therefore inequalities of this type enable one to tackle the discrete polynomials
similarly to the continuous ones. The classical approach (see e.g. \cite{ismail}, \cite{li}), 
heavily depends on the fact that the corresponding generating functions are of a rather simple form
allowing asymptotic investigation. 
As far as we know the method of Laguerre type inequalities is the only one
producing sharp explicit bounds in the discrete case. 
We outline this type of applications for the case of Krawtchouk polynomials in the last section.\\
The following observation
(a prototypical result in the case of Krawtchouk polynomials was given in
\cite{chih}) shows that a conception of mesh arises very naturally in this context.
\begin{theorem}
\label{mesh}
Let $p \in \RP$ satisfies (\ref{rec1}) and have all its roots in the
open interval $I$. Then $M(p) \ge 1$
provided  $c_k(x) > 0$ for $x \in I.$
If in addition $b_k(x) >0$ on $I ,$ then  $M(p) \ge 2.$
\end{theorem}
In view of this it would be very useful to have inequalities valid for $M(p) \ge 1.$
At present we know only how to "mend" (\ref{glform})
for $m=1,$ namely we prove
\begin{theorem}
\label{disin1}
Let $p \in$ \RP and $M(p) \ge 1$ then for any $\mu (x) \ge 0,$
\be
\label{disine}
U(p)=p^2(x)-p(x-1)p(x+1)+\frac{1}{4} \left( p(x+1)-\mu (x) p(x)+p(x-1) \right)^2 \ge 0
\ee
\end{theorem}
In principle one can apply (\ref{disine}) and an approach developed in \cite{fk1} for
obtaining bounds corresponding to the second term of asymptotics,
but it leads to rather cumbersome calculations in comparing with
(\ref{glform}). 
Yet (\ref{disine}) may be quite useful in some tough cases, e.g. for bounding
the zeros of Charlier polynomials.
The details are rather lengthy and will be given in a separate paper.
Finally, it worth noticing that some sharper versions
of the Laguerre inequality are known in the polynomial case. 
It has been shown  \cite{love}, \cite{fk3}
that if $p \in \RP$ has degree
at most $k$ then
$$(k-1)p'^2-k p p'' \ge 0,$$
$$3(k-2)(k-3)p''^2-4(k-1)(k-3)p'p'''+k(k-1)p p^{(4)} \ge 0,$$ 
It would be interesting to extend these results
on the discrete case. We suggest the following conjecture.
\begin{conjecture}
Let $p \in$ \RP, $deg(p)=k,$ and
$M(p) \ge 1,$ then
$$(k-1)(p(x+1)-p(x-1))^2-4k p(x)(p(x+1)-2p(x)+p(x-1)) \ge 0$$
\end{conjecture}

\section{Proofs} 
We split the proof of Theorem \ref{disin4} into several lemmas.
\begin{lemma}
\label{lem1}
If $p \in$ \RP then
\be
\label{delay}
V_m((x-a)p)= ((x-a)^2-m^2)V_m (p)+V_{m-1}(p)
\ee
\end{lemma}
\proof
$$V_m((x-a)p)=\sum_{j=-m}^m (-1)^j \frac{p(x-j)p(x+j)((x-a)^2-j^2)}{( m-j)!(m+j)!}=$$
$$\sum_{j=-m}^m (-1)^j \frac{p(x-j)p(x+j)((x-a)^2-m^2+m^2-j^2)}{( m-j)!(m+j)!}=$$
$$((x-a)^2-m^2)V_m (p)+\sum_{j=-m}^m (-1)^j \frac{p(x-j)p(x+j)(m^2-j^2)}{( m-j)!(m+j)!}=$$
$$((x-a)^2-m^2)V_m (p)+ V_{m-1}(p)$$
\QED

\begin{lemma}
\label{lem2}
Let $p=\prod_{i=1}^k (x-a_i)$, then $V_m (p)=0$ for $k<m$, $V_m(p)=1$ for $k=m$, and
\be
\label{glin}
V_m(p)= (m+1)x^2-2x \sum_{i=1}^{m+1} a_i+\sum_{i=1}^{m+1} a_i^2-\frac{1}{4} {2m+2 \choose 3}
\ee
for $k=m+1.$
\end{lemma}
\proof
For $m=1$ it is checked directly.
The result easily follows from (\ref{delay}) by the induction on $m$.
\QED

$Proof$ $of$ $Theorem$ \ref{disin4}:
It is enough to prove the theorem for $p \in$\RP. The general
result follows by $V_m(f e^{ \beta x})=e^{ \beta x}V_m(f),$
and easy limiting arguments.
Let $p=\prod_{i=1}^k (x-a_i),$ $a_1 <a_2 <...<a_k,$ and 
$M(p) \ge \delta_m =\sqrt{4-\frac{6}{m+2}}$. 
By Lemma \ref{lem2} the claim is true for $k \le m.$ We apply the induction on $k$.
Consider two cases:\\
$Case \, 1:$ $k=m+1$\\
It is enough to show that (\ref{glin}), a quadratic in $x$,  has a non-positive discriminant,
that is
$$\Delta = 4 \left( \sum_{i=1}^{m+1} a_i \right)^2-4(m+1)\sum_{i=1}^{m+1} a_i^2 +(m+1){2m+2 \choose 3} \le 0 .$$
Since by the assumption $a_i-a_{i-1} \ge \delta_m ,$ $i=2,..,m+1,$ we may 
put $a_i=(i-1) \delta_m+\sum_{j=1}^i r_j,$ where all $r_j \ge 0.$
It is easy to check that $\frac{\partial \Delta}{\partial r_i} < 0$ for $ 2 \le i \le m+1$
and $\frac{\partial \Delta}{\partial r_1} = 0.$
Hence to find the maximum value of $\Delta$ we can take all $r_i=0,$
yielding $a_i =(i-1) \delta_m$ and
$\Delta= 0.$ \\
$Case \, 2:$ $k>m+1$\\
Fix any $x=x_0,$
then $\max (|x_0-a_1|,|x_0-a_k|) \ge (k-1) M(p)/2,$ and so
for $a$ equals either $a_1$ or $a_k$ we have
$$(x_0-a)^2-m^2 \ge \frac{(2m+1)(k-1)^2}{2m+4}-m^2 \ge \frac{m^2+4m+1}{2m+4} >0.$$
On putting $p(x)=(x-a)q(x)$ the result follows from (\ref{delay}) by the induction hypothesis.
\QED

To prove Theorem \ref{disin1} we need
a more detailed version of Theorem \ref{disin4} for $m=1.$
\begin{lemma}
\label{lemd:1}
If $M(p) \ge 1, \, deg(p)=k \ge 2,$ and for some $x=x_0,$
$$p^2 (x_0)-p(x_0+1)p(x_0-1) <0,$$
then there are two consecutive roots $x_i ,x_{i+1}$ of $p(x)$
s.t. $x_0 \in (x_i ,x_{i+1})$ and $x_{i+1}-x_i < \sqrt{2}.$
\end{lemma}
\proof
Put $p(x)=a \prod_{j=1}^k (x-x_j), \; \; x_1 <...<x_k.$ If $p^2 (x_0)-p(x_0+1)p(x_0-1) <0$ then
$$1 < \frac{p(x_0+1)p(x_0-1)}{p^2(x_0)}=\prod_{j=1}^k (1-\frac{1}{(x_0-x_j)^2}):=
\prod_{j=1}^k q_j.$$
Since $|q_j|>1$ implies $q_j <-1$ and $|x_0-x_j|< \frac{1}{\sqrt{2}},$ 
it follows by $M(p) \ge 1$ and $\prod_{j=1}^k q_j >0,$ that $x_i <x_0 < x_{i+1}$
for some $i$ and $q_i q_{i+1} >1.$
But elementary calculations show that 
$q_i q_{i+1}=(1-\frac{1}{(x_0-x_i)^2})(1-\frac{1}{(x_0-x_{i+1})^2})>1$ is possible
only for $x_{i+1}-x_i \le \sqrt{2}.$
\QED
 
$Proof$ $of$ $Theorem$ \ref{disin1}:
To prove (\ref{disine}) assume $U(p) <0$ for some $x=x_0$,
Then $p(x+1)$ and $p(x-1)$ are of the same sign. Then by the
previous lemma $p(x)$ is of the opposite sign.
But on expanding $V_2$ we convince that then $U \ge 0.$
\QED

$Proof$ $of$ $Theorem$ \ref{mesh}:
Assume that $y_1 <y_2$ is the largest pair of consecutive roots of $p=p_k$ such that $y_2-y_1 <1.$
Since it must be a spectral point between any two roots of $p$ then
there is unique integer $m$, $y_1 \le m \le y_2$. Assume first that $m \ne y_1, \, y_2,$
then
$$Sign{(p(y_2+1))}=Sign(p(m+1))=-Sign(p(m)),$$
and
$$Sign(p(y_2-1))=Sign(p(m-1))=-Sign(p(m)),$$
contradicting $c_k(y_2) >0$ and $p(y_2+1) =-c_k(y_2)p(y_2-1).$
Cases $m=y_1$ or $y_2$ are similar.
Thus, we get $M(p) \ge 1.$\\
Let now $b_k(x) >0$, and $y_1 <y_2$ be two consecutive roots with $y_2-y_1 <2$.
Let $z=(y_1+y_2)/2$, then $Sign(p(z-1))=-Sign(p(z))=Sign(p(z+1)),$
contradicting (\ref{rec1}).
\QED

\section{Some Applications}
Inequalities (\ref{glform})
can be used to establish very sharp explicit bounds on discrete
orthogonal polynomials. The case of $L_1$ and $V_1$ for Hermite and Krawtchouk
polynomials respectively has been considered
in \cite{fk1},\cite{k3}. In the continuous case $L_2$ leads to much more precise estimates
than $L_1$, \cite{fk2}, \cite{kbes}.
Here we will show that $V_2$ gives a similar improvement for
the case of binary Krawtchouk polynomials
$P_k (x)=P_k^n (x), \, \, k<n/2.$ 
We outline the method leaving aside the details of calculations which are similar to these in
\cite{fk1}, \cite{fk2},
and require a symbolic package, we used Mathematica.
Asymptotic results obtained via a more classical approach 
can be found in \cite{ismail}, \cite{li}.\\
Let $x_1 <...<x_k ,$ be the zeros of $P_k (x)$.
Notice that $P_k (x)$ and consequently its zeros are symmetric with respect to $n/2 .$
First, we need rough bounds on $x_1$ and $x_k.$
Since Krawtchouk polynomials are orthogonal on $[0,n]$, we get $0 \le x_1 <x_k \le n.$
Slightly more accurate considerations yield
$1 < x_1 < x_k < n-1,$
and this is the bound we start with.
Thus, in the sequel we assume $1 < x < n-1.$
A difference equation for Krawtchouk polynomials is
\be
\label{krawrec}
(n-x)P_k(x+1)=(n-2k)P_k (x)-x P_k (x-1)
\ee
Since in our case $M(P_k(x)) >2$ by the assumption $k<n/2 ,$  one gets
$$V_2(P_k(x))= \frac{A(x) t^2+B(x) t +C(x)}{12(n-x)(n-x-1)(x-1)} (P_k(x))^2 \ge 0,$$
where $t=t(x)=P_k(x-1)/P_k(x),$ 
$$A(x)= - x(4x^2-4n x+4n+m^2-4)$$
$$B(x)=m(4x^2-4n x+2x+3n+m^2-4)$$
$$C(x)=4x^3-8n x^2+(4n^2+2n+m^2-4)x-2n^2-m^2 n+4n-m^2$$
and $m=n-2k.$
Since $t(x)$ attains all values from $- \infty$ to $\infty$ between any two consecutive zeros of
$P_k$ in the oscillatory region one must have 
$A(x)>0.$
It worth noticing that $A(x_1)>0$ gives
$x_1 > \frac{n}{2}-\sqrt{(k-1)(n-k-1)},$ in fact the first term of the corresponding asymptotic.
Thus, for $\Delta (x)= B^2(x)-4A(x)C(x) \ge 0,$ 
we get that $t(x)$ is confined between two solutions of
$A(x) t^2+B(x) t +C(x)=0.$
The explicit result is rather messy, therefore we only notice that
in fact the smaller branch yields very sharp approximation to $P_k(x-1)/P_k(x).$\\
To obtain a stronger inequality for the extreme zeros
observe that the function $t(x)$ is monotonically increases 
on $(1,x_1)$ from $P_k(0)/P_k(1)$
to $\infty $
and avoids the region given by 
$A(x) t^2+B(x) t +C(x)<0. $ This implies that 
the discriminant of this quadratic in $t$
is negative in the oscillatory region, i.e.
$\Delta (x_1) < 0.$
We have
\be
\label{krawr4}
\Delta(x)=(y^2-(n-1)^2+m^2-1)^3-2(y^2+m^2-1)^2+m^2 y^2+2(n-1)^2(n^2-2n+5) <0,
\ee
where $y=n-2x.$
The discriminant surface of this cubic equation in $z=y^2$ (that is the domain of parameters
where the equation has a multiple root)
is given by
$$(256(n-1)^2-27m^4+32m^2)(n^2-m^2)((n-2)^2-m^2)=0. $$
It easily implies that for $2 \le k < \frac{n}{2}-2 \cdot 3^{-3/4} \sqrt{n},$
equation (\ref{krawr4}) has two real roots, the positive one satisfying the inequality
$$y < 2 \sqrt{k(n-k)} \, \left( 1-\frac{1}{4} \left( \frac{n-2k}{k(n-k)} \right)^{2/3} \right)$$
Thus we conclude that for
$2 \le k < \frac{n}{2}-2 \cdot 3^{-3/4} \sqrt{n}$
all the roots of $P_k^n (x)$ are in the interval
\be
\label{rootk}
\frac{n}{2 } \pm
\sqrt{k(n-k)} \, \left( 1-\frac{1}{4} \left( \frac{n-2k}{k(n-k)} \right)^{2/3} \right)
\ee
Slightly more accurate bounds
$\frac{3}{8}$ instead of $\frac{1}{4}$, but only for $k$ growing linearly with
$n$, has been found in \cite{fk1}. It seems that the exact (asymptotic) value of the
constant at the error term is unknown, although it almost certainly relates to
the Airy function.\\
Now we will show how to estimate $V_2$ in the oscillatory region.
Bounds on $P_k (x)$ and $P_k (x-1)$ inside the oscillatory region readily follow from the
corresponding bounds on $V_2 ,$ as the major axes of the ellipse corresponding to $V_2.$
The situation in the transitory regions, that is
near to the extreme
zeros, is more complicated. 
It either requires an appealing to Remez inequality or
some more tricky arguments, see e.g.  \cite{fk3, k3}.\\
The expression
\be
\label{wform}
W(x)=\frac{V_2(P_k(x+1))-z V(P_k(x))}{(P_k(x))^2}
\ee
is a quadratic in $t$, more precisely it can be written as
$$\frac{C(x)t^2+D(x)t+E(x)}{(n-x-2)(n-x-1)(n-x)^2}, $$
where $C(x),D(x),E(x)$ are certain polynomials.
Choose $z$ from the condition $D^2(x)-4C(x)E(x)=0,$ 
then the sign of $W(x)$ coincides with the sign of $C(x)$.
Calculations yield
$$z=\frac{(x-1)(\Delta (x+\frac{1}{2} )-3S \pm
6 \sqrt{R} \, )}{(n-x-2) \Delta (x)}$$
where 
$$S=y(y-2)(y-1)^2-(n^2-2n-m^2+2)^2+7(n-1)^2,$$
$$R=(n^2-y^2-2n+2y)(n^2-m^2)((n-2)^2-m^2)((n-1)^2-m^2-(y-1)^2).$$
Thus, by $V_2 > 0,$ we obtain that for $2 \le k < \frac{n}{2}-2 \cdot 3^{-3/4} \sqrt{n},$
$$
\frac{(x-1)(\Delta (x+\frac{1}{2} )-3S +
6 \sqrt{R} \, )}{(n-x-2) \Delta (x)} \le \frac{ V_2 (P_k(x+1))}{V_2 (P_k(x))} \le
\frac{(x-1)(\Delta (x+\frac{1}{2} )-3S -
6 \sqrt{R} \, )}{(n-x-2) \Delta (x) } ,
$$
provided $\Delta (x)<0.$\\
In our case $V_2 (P_k ( \frac{n}{2}))$ can be found explicitly, giving the required initial conditions.
It is important to notice that for any fixed $-1 < \epsilon <1,$
and $x=\frac{n}{2}+2 \epsilon \sqrt{k(n-k)},$ we have
$\sqrt{R}/\Delta(x)=O(n^{-2}).$
In fact this can be improved to $O(n^{-3})$ 
if we consider 
$$V_2(P_k(x))+\frac{(P_k(x+1)-P_k(x-1))^2}{2(4n x-(n-2k)^2-4x^2)}$$
in (\ref{wform}) instead of $V_2(P_k(x))$, that is
perturbing the form in
a similar way to that used in \cite{fk2} for the case of Hermite polynomials.

%

\end{document}